\documentclass[matprg]{mcsreport}
\usepackage{subeqn}
\usepackage{extra}

\newcommand{\tlam}{\tilde{\lambda}}
\newcommand{\hlam}{\hat{\lambda}}

\newcommand{\nabg}{\nabla g}

\newcommand{\cSlam}{\cS_{\lambda}}
\newcommand{\epslam}{\epsilon_{\lambda}}
\newcommand{\zsl}{z_{[\ell]}}
\newcommand{\lsl}{\lambda_{[\ell]}}

\newcommand{\comment}[1]{}

\title{Constraint Identification and Algorithm Stabilization 
  for Degenerate Nonlinear Programs}

\author{Stephen J. Wright%
\thanks{Research supported by the Mathematical,
Information, and Computational Sciences Division subprogram of the
Office of Advanced Scientific Computing Research, U.S. Department of
Energy, under Contract W-31-109-Eng-38.}}

\institute{Stephen J. Wright\at Mathematics and Computer Science Division,
Argonne National Laboratory, 9700 South Cass Avenue, 
Argonne, Illinois 60439; {\tt wright@mcs.anl.gov}}

\date{\today}
\reportnumber{P865-1200, December, 2000}

\titlerunning{Constraint Identification for Degenerate Nonlinear
Programs}

\begin{document}

\maketitle

\begin{abstract}
  In the vicinity of a solution of a nonlinear programming problem at
  which both strict complementarity and linear independence of the
  active constraints may fail to hold, we describe a technique for
  distinguishing weakly active from strongly active constraints.  We
  show that this information can be used to modify the sequential
  quadratic programming algorithm so that it exhibits
  superlinear convergence to the solution under assumptions weaker
  than those made in previous analyses. 
%
%
\end{abstract}



\section{Introduction} \label{introduction}

Consider the following nonlinear programming problem with inequality
constraints:
\beq \label{nlp}
\makebox{\rm NLP:} \hspace*{.7in} 
\min_z \, \phi(z) \gap \makebox{\rm subject to $g(z) \le 0$},
\eeq
where $\phi:\R^n \to \R$ and $g:\R^n \to \R^m$ are twice Lipschitz
continuously differentiable functions.  Optimality conditions for
\eqnok{nlp} can be derived from the
Lagrangian for \eqnok{nlp}, which is
\beq \label{lagr}
\cL(z,\lambda) = 
\phi(z) + \lambda^T g(z),
\eeq
where $\lambda \in \R^m$ is the vector of Lagrange multipliers.  When
a constraint qualification holds at $z^*$ (see discussion below), the
first-order necessary conditions for $z^*$ to be a local solution of
\eqnok{nlp} are that there exists a vector $\lambda^* \in \R^m$ such
that
\beq \label{kkt}
\cL_z (z^*,\lambda^*) = 0, \sgap g(z^*) \le 0, \sgap  \lambda^* \ge 0, \sgap
(\lambda^*)^T g(z^*) = 0.
\eeq
These relations are the well-known Karush-Kuhn-Tucker (KKT)
conditions. The set $\cB$ of active constraints at $z^*$ is 
\beq \label{act.def}
\cB = \{ i=1,2,\dots,m \, | \, g_i(z^*) = 0 \}.
\eeq
It follows immediately from \eqnok{kkt} that we can have $\lambda^*_i
>0$ only if $i \in \cB$. The {\em weakly active} constraints are identified
by the indices $i \in \cB$ for which $\lambda^*_i=0$ for all
$\lambda^*$ satisfying \eqnok{kkt}. Conversely, the {\em strongly active}
constraints are those for which $\lambda^*_i>0$ for at least one
multiplier $\lambda^*$ satisfying \eqnok{kkt}. The strict
complementarity condition holds at $z^*$ if there are no weakly active
constraints.

We are interested in degenerate problems, those for which the active
constraint gradients at the solution is linearly dependent or the
strict complementarity condition fails to hold (or both).  The first
part of our paper describes a technique for partitioning $\cB$ into
weakly active and strongly active indices. 
Section~\ref{sec:detect} builds on the technique described by
Facchinei, Fischer, and Kanzow~\cite{FacFK98} for identifying $\cB$. Our
technique requires the solution of a sequence of closely related
linear programming subproblems in which the set of strongly active
indices is assembled progressively. Solution of one additional linear
program yields a Lagrange multiplier estimate $\lambda$ such that the
components $\lambda_i$ for all strongly active indices $i$ are bounded
below by a positive constant.

In the second part of the paper, we use the cited technique to adjust
the Lagrange multiplier estimate between iterations of the stabilized
sequential quadratic programming (sSQP) algorithm described by
Wright~\cite{Wri98a} and Hager~\cite{Hag97a}. The resulting technique
has the advantage that it converges superlinearly under weaker
conditions than considered in these earlier papers. We can drop the
assumption of strict complementarity and a ``sufficiently interior''
starting point made in \cite{Wri98a}, and we do not need the stronger
second-order conditions of \cite{Hag97a}.  Motivation for the sSQP
approach came from work on primal-dual interior-point algorithms
described in \cite{WriR94,RalW96b}. It is also closely related to the
method of multipliers and the ``recursive successive quadratic
programming'' approach of Bartholomew-Biggs~\cite{Big87}. (See
Wright~\cite[Section~6]{Wri97e} for a discussion of the similarities.)

Other work on stabilization of the SQP approach to yield superlinear
convergence under weakened conditions has been performed by
Fischer~\cite{Fis99a} and Wright~\cite{Wri97e}. Fischer
proposed an algorithm in which an additional quadratic program is
solved between iterations of SQP in order to adjust the Lagrange
multiplier estimate. He proved superlinear convergence under
conditions that are weaker than the standard assumptions but stronger
than the ones made in this paper.  Wright described
superlinear local convergence properties of a class of inexact SQP
methods and showed that sSQP and Fischer's method could be expressed
as members of this class. This paper also introduced a modification of
standard SQP that enforced only a subset of the linearized
constraints---those in a ``strictly active working set''---and
permitted slight violations of the nonenforced constraints yet
achieved superlinear convergence under weaker-than-usual conditions.

Bonnans~\cite{Bon89a} showed that when strict complementarity fails to
hold but the active constraint gradients are linearly independent,
then the standard SQP algorithm (in which any nonuniqueness in the
solution of the SQP subproblem is resolved by taking the solution of
minimum norm) converges superlinearly.

Our concern here is with {\em local} behavior, so we assume
availability of a starting point $(z^0, \lambda^0)$ that is
``sufficiently close'' to the optimal primal-dual set. We believe,
however, that ingredients of the approach proposed here can be
embedded in practical algorithms, such as SQP algorithms that include
modifications (merit functions and filters) to ensure global
convergence. We believe also that this approach could be used to
enhance the robustness and convergence rate of other types of
algorithms, including augmented Lagrangian and interior-point
algorithms, in problems in which there is degeneracy at the
solution. We mention one such extension in Section~\ref{sec:ext}.

\section{Assumptions, Notation, and Basic Results} \label{assumptions}

We now review the optimality conditions for \eqnok{nlp} and outline
the assumptions that are used in subsequent sections. These include
the second-order sufficient condition we use here, the
Mangasarian-Fromovitz constraint qualification, and the
definition of weakly-active indices.

Recall the KKT conditions \eqnok{kkt}. The set of ``optimal'' Lagrange
multipliers $\lambda^*$ is denoted by $\cSlam$, and the primal-dual
optimal set is denoted by $\cS$. Specifically, we have
\beq \label{def.SS}
\cSlam = \{ \lambda^* \, | \, \makebox{\rm $\lambda^*$ satisfies
(\protect\ref{kkt})} \},  \sgap
\cS = \{ z^* \} \times \cSlam.
\eeq
An alternative, compact form of the KKT conditions is the following 
variational inequality formulation:
\beq \label{kkt.N}
\bmat{c} \nabla \phi(z^*) + \nabg(z^*) \lambda^* \\ g(z^*) \emat \in 
\bmat{c} 0 \\ N(\lambda^*) \emat,
\eeq
where $N(\lambda)$ is the set defined by
\beq \label{def.Nl}
N(\lambda) \defeq \left\{
\begin{array}{cl} \{ y \, | \, y \le 0 \; 
\mbox{\rm and} \; y^T \lambda = 0 \} & \makebox{\rm if $\lambda \ge 0$,} \\
\emptyset & \makebox{\rm otherwise.}
\end{array}
\right.
\eeq

We now introduce notation for subsets of the set $\cB$ of active
constraint indices at $z^*$, defined in \eqnok{act.def}.  For any
optimal multiplier $\lambda^* \in \cSlam$, we define the set
$\cB_+(\lambda^*)$ to be the ``support'' of $\lambda^*$, that is,
\[
\cB_+(\lambda^*) = \{  i \in \cB \, | \, \ \lambda^*_i>0 \}.
\]
We define $\cB_+$ (without argument) as
\beq \label{def.Bp}
\cB_+ \defeq \cup_{\lambda^* \in \cSlam} \, \cB_+ (\lambda^*);
\eeq
this set contains the indices of the {\em strongly active}
constraints. Its complement in  $\cB$ is denoted by $\cB_0$, that is,
\[
\cB_0 \defeq \cB \backslash \cB_+.
\]
This set $\cB_0$ contains the {\em weakly active} constraint indices,
those indices $i \in \cB$ such that $\lambda^*_i=0$ for all $\lambda^*
\in \cSlam$.  In later sections, we make use of the quantity $\epslam$
defined by
\beq \label{def.dellam}
\epslam \defeq \max_{\lambda^* \in \cSlam}  \min_{i \in \cB_+} 
\lambda^*_i.
\eeq
Note by the definition of $\cB_+$ that $\epslam>0$.


The Mangasarian-Fromovitz constraint qualification
(MFCQ)~\cite{ManF67} holds at $z^*$ if there is a vector $\bar{y} \in
\R^{n}$ such that
\[
\nabg_i (z^*)^T \bar{y} <0 \sgap \makebox{\rm for all $i \in \cB$.}
\]
By defining $\nabg_{\cB}$ to be the $n \times | \cB|$ matrix whose
rows are $\nabla g_i(\cdot)$, $i \in \cB$, we can write this condition
alternatively as
\beq \label{mfcq}
\nabg_{\cB}(z^*)^T \bar{y} <0.
\eeq
It is well known that MFCQ is equivalent to boundedness of the set
$\cSlam$; see Gauvin~\cite{Gau77}.

Since $\cSlam$ is defined by the linear conditions $\nabla \phi(z^*) +
\nabg(z^*) \lambda^*$ and $\lambda^* \ge 0$, it is closed and convex.
Therefore, under MFCQ, it is also compact.

We assume throughout that the following second-order condition is
satisfied: there is $\sigma>0$ such that
\beq \label{2os}
w^T \cL_{zz}(z^*,\lambda^*) w  \ge \sigma \| w \|^2, \sgap
\makebox{\rm for all  $\lambda^* \in \cSlam$},
\eeq
and for all $w$ such that
\beq \label{2os.w}
\begin{array}{ll} 
\nabg_i(z^*)^T w =0, & \makebox{\rm for all $i \in \cB_+$}, \\
\nabg_i(z^*)^T w \le 0, & \makebox{\rm for all $i \in \cB_0$}.
\end{array}
\eeq
This condition is referred to as Condition 2s.1 in
\cite[Section~3]{Wri97e}. 
Weaker second-order conditions, stated in terms of a quadratic growth
condition of the objective $\phi(z)$ in a feasible neighborhood of
$z^*$, are discussed by Bonnans and Ioffe~\cite{BonI95} and
Anitescu~\cite{Ani99b}.

Our standing assumption for this paper is as follows.
\begin{assumption} \label{ass:standing}
  The first-order conditions \eqnok{kkt}, the MFCQ \eqnok{mfcq}, and
  the second-order condition \eqnok{2os}, \eqnok{2os.w} are satisfied
  at $z^*$. Moreover, the functions $\phi$ and $g$ are twice Lipschitz
  continuously differentiable in a neighborhood of $z^*$.
\end{assumption}

\noindent
The following is an immediate consequence of this assumption.
\begin{theorem} \label{th:2os}
  Suppose that Assumption~\ref{ass:standing} holds. Then $z^*$ is an
  isolated stationary point and a strict local minimizer of
  \eqnok{nlp}.
\end{theorem}
\begin{proof}
See Robinson~\cite[Theorems~2.2 and 2.4]{Rob82}.
\end{proof}

We use the notation $\delta(\cdot)$ to denote  distances from
the primal, dual, and primal-dual optimal sets, according to
context. Specifically, we define
\beq \label{errdef}
\delta(z) \defeq \| z-z^* \|, \sgap
\delta(\lambda) \defeq \dist (\lambda, \cSlam), \sgap
\delta(z,\lambda) \defeq \dist ((z,\lambda), \cS),
\eeq
where $\| \cdot \|$ denotes the Euclidean norm unless a subscript
specifically indicates otherwise.
We also use $P(\lambda)$ to denote the projection of $\lambda$ onto
$\cSlam$; that is, we have $P(\lambda) \in \cSlam$ and
$\| P(\lambda) - \lambda \| = \dist (\lambda, \cSlam)$.  Note
that from \eqnok{errdef}  we have $\delta(z)^2 + \delta(\lambda)^2 =
\delta(z,\lambda)^2$, and therefore 
\beq \label{deltas}
\delta(z) \le \delta(z,\lambda), \gap
\delta(\lambda) \le \delta(z,\lambda).
\eeq

Using Assumption~\ref{ass:standing}, we can prove the following
result, which gives a practical way to estimate the distance
$\delta(z,\lambda)$ of $(z,\lambda)$ to the primal-dual solution set
$\cS$.
\begin{theorem} \label{th:dlamz}
  Suppose that Assumption~\ref{ass:standing} holds. Then there are
  positive constants $\delta$, $\kappa_0$, and $\kappa_1$ such that
  for all $(z,\lambda)$ with $\delta(z,\lambda) \le \delta$, the
  quantity $\eta(z,\lambda)$ defined by
\beq \label{def.eta}
\eta(z,\lambda) \defeq \left\| \bmat{c} \cL_z(z,\lambda)  \\
\min(\lambda, -g(z)) \emat \right\|
\eeq
(where $\min(\lambda, -g(z))$ denotes the vector whose
$i$th component is $\min(\lambda_i, -g_i(z))$) satisfies
\[
\kappa_0 \delta(z,\lambda) \le \eta(z,\lambda) \le \kappa_1 \delta(z,\lambda).
\]
\end{theorem}

\noindent
See Facchinei, Fischer, and Kanzow~\cite[Theorem~3.6]{FacFK98},
Wright~\cite[Theorem~A.1]{Wri97e}, and Hager and
Gowda~\cite[Lemma~2]{HagG97} for proofs of this result. (The
second-order condition is stated in a slightly different fashion in
\cite{FacFK98} but is equivalent to \eqnok{2os}, \eqnok{2os.w}.)

We use order notation in the following (fairly standard) way: If two matrix,
vector, or scalar quantities $M$  and $A$ are
functions of a common quantity,  we write $M=O(\|A \|)$ if there is a
constant $\beta$ such that $\| M \| \le \beta \|A \|$ whenever $\|A \|$ is 
sufficiently small. We write $M=\Omega (\|A \|)$ if there is a
constant $\beta$ such that $\| M \| \ge \beta^{-1} \|A \|$ whenever
$\|A \|$ sufficiently small, and $M = \Theta (\|A \|)$ if both
$M=O(\|A \|)$ and $M=\Omega (\|A \|)$. We write $M = o(\| A \|)$ if
for all sequences $\{A_k\}$ with $\| A_k \| \to 0$, the corresponding
sequence $\{M_k\}$ satisfies $\| M_k \| / \| A_k \| \to 0$.
By using this notation, we can rewrite the conclusion of
Theorem~\ref{th:dlamz} as follows:
\beq \label{delest}
\eta(z,\lambda) = \Theta (\delta(z,\lambda)).
\eeq

\section{Detecting Active Constraints} \label{sec:detect}

We now describe a procedure, named Procedure ID0, for identifying
those inequality constraints that are active and the solution, and
classifying them according to whether they are weakly active or
strongly active. We prove that Procedure ID0 classifies the indices
correctly given a point $(z,\lambda)$ sufficiently close to the
primal-dual optimal set $\cS$. Finally, we describe some
implementation issues for this procedure.

\subsection{The Detection Procedure} \label{subsec:proc}

Facchinei, Fischer, and Kanzow~\cite{FacFK98} showed that the function
$\eta(z,\lambda)$ defined in \eqnok{delest} can be used as the basis
of a scheme for identifying the active set $\cB$.  Choosing some
$\tau \in (0,1)$, they estimated 
\beq \label{def.Best}
\cA(z,\lambda) \defeq 
\{ i=1,2,\dots,m \, | \, g_i(z) \ge -\eta(z,\lambda)^{\tau} \}.
\eeq
We have the following result.
\begin{theorem} \label{th:Best}
  Suppose that Assumption~\ref{ass:standing} holds. Then there exists
  $\delta>0$ such that for all $(z,\lambda)$ with $\delta(z,\lambda)
  \le \delta$, we have $\cA(z,\lambda) = \cB$.
\end{theorem}
\begin{proof}
  The result follows immediately from \cite[Definition~2.1,
  Theorem~2.3]{FacFK98} and Theorem~\ref{th:dlamz} above.
\end{proof}

A scheme for estimating $\cB_+$ (hence, $\cB_0$) is described in
\cite{FacFK98}, but it requires the strict MFCQ condition to hold,
which implies that $\cSlam$ is a singleton. Here we describe a more
complicated scheme for estimating $\cB_+$ that requires only the
conditions of Theorem~\ref{th:Best} to hold.

Our scheme is based on linear programming subproblems of the following
form, for a given parameter ${\tau} \in (0,1)$ and a given set $\hat{\cA}
\subset \cA(z,\lambda)$:
\begin{subequations} \label{B0est}
\beqa \label{B0est.1}
& \max_{\tlam} \, 
\sum_{i \in {\hat{\cA}}} \tlam_i \;\; \mbox{subject to} \\
\label{B0est.2}
& -\eta(z,\lambda)^{\tau} \le 
\nabla \phi(z) + \sum_{i \in \cA(z,\lambda)} \tlam_i \nabla g_i(z)
\le \eta(z,\lambda)^{\tau} \\
\label{B0est.3}
& \tlam_i  \ge 0, \;\; \mbox{for all} \; i \in \cA(z,\lambda); \;\;\;
\tlam_i=0 \;\; \mbox{otherwise}.
\eeqa
\end{subequations}
Note that the objective function involves elements $\tlam_i$ only for
indices $i$ in the subset $\hat{\cA}$, whereas the $\tlam_i$ are
permitted to be nonzero for all $i \in \cA(z, \lambda)$. The idea is
that $\hat{\cA}$ contains those indices that {\em may} belong to
$\cB_0$; by the time we solve \eqnok{B0est}, we have already decided
that the other indices $i \in \cA(z,\lambda) \backslash \hat{\cA}$
probably belong to $\cB_+$.

The complete procedure is as follows.

\medskip

\btab
\> {\bf Procedure ID0} \\
\> Given constants $\tau$ and $\hat{\tau}$ 
satisfying $0 < \hat{\tau} < \tau < 1$, and point $(z,\lambda)$; \\
\> Evaluate $\eta(z,\lambda)$ from \eqnok{def.eta} and $\cA(z,\lambda)$ from 
\eqnok{def.Best}; \\
\> Define $\hat{\cA}_{\rm init} = \cA(z,\lambda) \backslash \{ i \, | \, 
\lambda_i \ge \eta(z,\lambda)^{\hat{\tau}} \}$; \\
\> $\hat{\cA} \leftarrow \hat{\cA}_{\rm init}$; \\
\> {\bf repeat} \\
\>\> solve \eqnok{B0est} to find $\tlam$; \\
\>\> set $\cC = \{ i \in \hat{\cA} \, | \, \tlam_i \ge \eta(z,\lambda)^{\hat{\tau}} \}$; \\
\>\> {\bf if} $\cC = \emptyset$ \\
\>\>\> stop with $\cA_0 = \hat{\cA}$, 
$\cA_+ = \cA(z,\lambda) \backslash \hat{\cA}$; \\
\>\> {\bf else} \\
\>\>\> set $\hat{\cA} \leftarrow \hat{\cA} \backslash \cC$; \\
\>\>\> {\bf if} $\hat{\cA} = \emptyset$ \\
\>\>\>\> stop with $\cA_0 = \emptyset$, $\cA_+ = \cA(z,\lambda)$; \\
\>\>\> {\bf end(if)} \\
\>\> {\bf end(if)} \\
\> {\bf end(repeat)}\\
\etab

\medskip

This procedure terminates finitely; in fact, the number of times the
``repeat'' loop executes is bounded by the cardinality of
$\hat{\cA}_{\rm init}$.

We prove that Procedure ID0 successfully identifies $\cB_+$ (for all
$\delta(z,\lambda)$ sufficiently small) in several steps, culminating
in Theorem~\ref{th:ID0.3}.
First, we estimate  the distance of $(z,\tlam)$ to the
solution set $\cS$, where $\tlam$ is the solution of \eqnok{B0est} for
some $\hat{\cA}$.
\begin{lemma} \label{lem:ID0.0}
  Suppose that Assumption~\ref{ass:standing} holds. Then there are
  positive constants $\delta_0$ and $\kappa_2$ such that whenever
  $\delta(z,\lambda) \le \delta_0$, any feasible point $\tlam$ of
  \eqnok{B0est} at any iteration of Procedure ID0 satisfies
\[
\delta(z,\tlam) \le \kappa_2 \delta(z,\lambda)^{\tau}.
\]
\end{lemma}
\begin{proof}
  Initially choose $\delta_0 = \delta$ for $\delta$ defined in
  Theorem~\ref{th:Best}, so that $\cA(z,\lambda) = \cB$. Hence, we
  have $\hat{\cA} \subset \cB$ at all iterations of Procedure ID0.
  
  We now estimate $\eta(z,\tlam)$ using the definition
  \eqnok{def.eta}. We have directly from the constraints
  \eqnok{B0est.2} that
\[
\| \cL_z(z,\tlam) \|_{\infty} \le \eta(z,\lambda)^{\tau}.
\]
For the vector $\min(\tlam, -g(z))$, we have for $i \in \cB$ that
$g_i(z^*)=0$ and $\tlam_i \ge 0$, and so
\[
i \in \cB \; \Rightarrow \; 
| \min(\tlam_i, -g_i(z)) | 
\le | g_i(z)| = O(\| z-z^*\|) = O(\delta(z,\lambda)).
\]
Meanwhile for $i \notin \cB= \cA(z,\lambda)$, we have $\tlam_i=0$ and
$g_i(z^*)<0$, and so
\[
i \notin \cB \; \Rightarrow \; 
| \min(\tlam_i, -g_i(z)) | = \max(0,g_i(z)) \le
|g_i(z) - g_i(z^*) | = O(\delta(z,\lambda)).
\]
By substituting these estimates into \eqnok{def.eta}, and using the
equivalence of $\| \cdot \|_{\infty}$ and the Euclidean norm and the
result of Theorem~\ref{th:dlamz}, we have that there is a constant
$\bar{\kappa}_2>0$ such that
\[
\eta(z,\tlam) \le \bar{\kappa}_2 \delta(z,\lambda)^{\tau}.
\]
Using Theorem~\ref{th:dlamz} again, we have
\beq \label{id0.1.5}
\delta(z,\tlam) \le \kappa_0^{-1} \eta(z,\tlam) \le
\kappa_0^{-1} \bar{\kappa}_2 \delta(z,\lambda)^{\tau},
\eeq
giving the result.
\end{proof}

In the next two lemmas and Theorem~\ref{th:ID0.3}, we show that for
$\delta(z,\lambda)$ sufficiently small, Procedure ID0 terminates with
$\cA_0 = \cB_0$ and $\cA_+ = \cB_+$.
\begin{lemma} \label{lem:ID0.1}
  Suppose that Assumption~\ref{ass:standing} holds. Then there is
  $\delta_1>0$ such that whenever $\delta(z,\lambda) \le \delta_1$,
  Procedure ID0 terminates with $\cB_0 \subset \cA_0$.
\end{lemma}
\begin{proof}
  Since we know the procedure terminates finitely, we need show only
  that $\cB_0 \subset \hat{\cA}$ at all iterations of the procedure.
  Initially set $\delta_1 = \delta_0 \le \delta$, so that
  $\cA(z,\lambda) = \cB$ and the result of Lemma~\ref{lem:ID0.0}
  holds. Suppose for contradiction there is an index $j \in \cB_0$
  such that $j$ either is not included in the initial index set
  $\hat{\cA}_{\rm init}$ or else is deleted from $\hat{\cA}$ at some
  iteration of Procedure ID0.
  
  Suppose first that $j$ is not included in $\hat{\cA}_{\rm init}$.
  Then we must have $\lambda_j > \eta(z,\lambda)^{\hat{\tau}}$, which
  by Theorem~\ref{th:dlamz} implies that 
\beq \label{id0.1.1}
\delta(z,\lambda) \ge | \lambda_j | \ge \eta(z,\lambda)^{\hat{\tau}} \ge
\kappa_0^{\hat{\tau}} \delta(z,\lambda)^{\hat{\tau}}.
\eeq
However, by decreasing $\delta_1$ and using $\hat{\tau} \in (0,1)$,
we can ensure that \eqnok{id0.1.1} does not hold whenever $\delta(z,\lambda)
\le \delta_1$. Hence, $j$ is included in $\hat{\cA}_{\rm init}$.

Suppose now that $j \in \cB_0$ is deleted from $\hat{\cA}$ at
some subsequent iteration. For this to happen, the subproblem
\eqnok{B0est} must have a solution $\tlam$ with 
\beq \label{id0.1.2}
\tlam_j > \eta(z,\lambda)^{\hat{\tau}}
\eeq
for some $\hat{\cA} \subset \cB$. Hence from Theorem~\ref{th:dlamz},
we have that
\beq \label{id0.1.4}
\delta(z,\tlam) \ge \tlam_j > \eta(z,\lambda)^{\hat{\tau}}
\ge \kappa_0^{\hat{\tau}} \delta(z,\lambda)^{\hat{\tau}}.
\eeq
By combining the result of Lemma~\ref{lem:ID0.0} with \eqnok{id0.1.4},
we have that
\[
\kappa_2 \delta(z,\lambda)^{\tau} \ge 
\kappa_0^{\hat{\tau}} \delta(z,\lambda)^{\hat{\tau}}.
\]
However, this inequality cannot hold when $\delta(z,\lambda)$ is
smaller than $(\kappa_0^{\hat{\tau}}
\kappa_2^{-1})^{1/(\tau-\hat{\tau})}$. Therefore, by decreasing
$\delta_1$ if necessary, we have a contradiction in this case also.
\end{proof}

\begin{lemma} \label{lem:ID0.2}
  Suppose that Assumption~\ref{ass:standing} holds. Then there is
  $\delta_2>0$ such that whenever $\delta(z,\lambda) \le \delta_2$,
  Procedure ID0 terminates with $\cB_+ \subset \cA_+$.
\end{lemma}
\begin{proof}
  Given any $j \in \cB_+$, we have for sufficiently small choice of
  $\delta_2$ that $j \in \cA(z,\lambda)$. We prove the result by
  showing that Procedure ID0 cannot terminate with $j \in \cA_0$.
  
  We initially set $\delta_2 = \delta_1$, where $\delta_1$ is the
  constant from Lemma~\ref{lem:ID0.1}.  (We reduce it as necessary,
  but maintain $\delta_2>0$, in the course of the proof.)  For
  contradiction, assume that there is $j \in \cB_+$ such that $j \in
  \hat{\cA}$ at all iterations of Procedure ID0, including the
  iteration on which the procedure terminates and sets $\cA_0 =
  \hat{\cA}$. Recalling the definition \eqnok{def.dellam} of
  $\epslam$, we use compactness of $\cSlam$ to choose $\lambda^* \in
  \cSlam$ such that $\epslam = \min_{i \in \cB_+} \lambda^*_i$. In
  particular, we have
\[
\lambda^*_j \ge \epslam>0
\]
for our chosen index $j$.  We claim that, by reducing $\delta_2$ if
necessary, we can ensure that $\lambda^*$ is feasible for
\eqnok{B0est} whenever $\delta(z,\lambda) \le \delta_2$. Obviously,
since $\cA(z,\lambda) = \cB$ by Theorem~\ref{th:Best}, $\lambda^*$ is
feasible with respect to \eqnok{B0est.3}.  Since $\lambda^* \in
\cSlam$ and
\[
\|z-z^*\| \le \delta(z,\lambda) \le \kappa_0^{-1} \eta(z,\lambda), 
\]
we have
\beqa
\nonumber
\left\|  \nabla \phi(z) + \sum_{i=1}^m \lambda^*_i \nabla g_i(z) 
\right\|_{\infty}
&=& 
\left\|  \nabla \phi(z) -\nabla \phi(z^*) 
+ \sum_{i=1}^m \lambda^*_i (\nabla g_i(z) -  \nabla g_i(z^*)) 
\right\|_{\infty} \\
\label{ID0.2.2}
& \le & M \| z-z^* \| \le M \kappa_0^{-1}  \eta(z,\lambda),
\eeqa
for some constant $M$ that depends on the norms of
$\nabla^2\phi(\cdot)$ and $\nabla^2 g_i(\cdot)$, $i \in \cB_+$ in the
neighborhood of $z^*$ and on a bound on the set $\cSlam$ (which {\em is}
bounded, because of MFCQ). Since $\tau<1$ and since $\eta(z,\lambda) =
\Theta (\delta(z,\lambda))$, we can reduce $\delta_2$ if necessary to
ensure that
\[
 M \kappa_0^{-1}  \eta(z,\lambda) < \eta(z,\lambda)^{\tau}
\]
whenever $\delta(z,\lambda) \le \delta_2$, thereby ensuring that the
constraints \eqnok{B0est.2} are satisfied by $\lambda^*$.

Since $\lambda^*$ is feasible for \eqnok{B0est}, a lower bound on the
optimal objective is 
\[
\sum_{i \in \hat{\cA}} \lambda^*_i \ge \lambda^*_j  \ge \epslam.
\]
However, since Procedure ID0 terminates with $j \in \hat{\cA}$, we
must have that $\cC = \emptyset$ for the solution $\tlam$ of
\eqnok{B0est} with this particular choice of $\hat{\cA}$. But we can
have $\cC = \emptyset$ only if $\tlam_i <
\eta(z,\lambda)^{\hat{\tau}}$ for all $i \in \hat{\cA}$, which means
that the optimal objective is no greater than $m
\eta(z,\lambda)^{\hat{\tau}}$. But since $\eta(z,\lambda) = \Theta
(\delta(z,\lambda))$, we can reduce $\delta_2$ if necessary to ensure
that
\[
m \eta(z,\lambda)^{\hat{\tau}} < \epslam
\]
whenever $\delta(z,\lambda) \le \delta_2$.  This gives a
contradiction, so that $\cA_0$ (which is set by Procedure ID0 to the
final $\hat{\cA}$) can contain no indices $j \in \cB_+$. Since $\cB_+
\subset \cB = \cA(z,\lambda)$ whenever $\delta(z,\lambda) \le
\delta_2$, we must therefore have $\cB_+ \subset \cA_+$, as claimed.
\end{proof}

By using the quantity $\delta_2$ from Lemma~\ref{lem:ID0.2}, we
combine this result with Theorem~\ref{th:Best} and
Lemma~\ref{lem:ID0.1} to obtain the following theorem.
\begin{theorem} \label{th:ID0.3}
Suppose that Assumption~\ref{ass:standing} holds. Then there is
  $\delta_2>0$ such that whenever $\delta(z,\lambda) \le \delta_2$,
  Procedure ID0 terminates with $\cA_+ = \cB_+$ and $\cA_0 = \cB_0$.
\end{theorem}

\subsection{Scheme for Finding an Interior Multiplier Estimate} \label{subsec:interior}

We now describe a scheme for finding a vector $\hlam$ that is close to
$\cSlam$ but not too close to the relative boundary of this set. In
other words, the quantity $\min_{i \in \cB_+} \hlam_i$ is not too far
from its maximum achievable value $\epslam$.

We find $\hlam$ by solving a linear programming problem similar to
\eqnok{B0est} but containing an extra variable to represent $\min_{i
\in \cB_+} \hlam_i$. We state this problem as follows:
\begin{subequations} \label{intlam}
\beqa \label{intlam.a}
& \max_{\hat{t}, \hlam} \, \hat{t} \;\; \mbox{subject to} \\
& \label{intlam.b}
\hat{t} \le \hlam_i, \;\; \mbox{for all $i\in \cA_+$}, \\
\label{intlam.c}
& -\eta(z,\lambda)^{\tau} e \le 
\nabla \phi(z) + \sum_{i \in \cA_+} \hlam_i \nabla g_i(z)
\le \eta(z,\lambda)^{\tau} e \\
\label{intlam.4}
& \hlam_i  \ge 0, \;\; \mbox{for all} \; i \in \cA_+; \;\;\;
\hlam_i=0 \;\; \mbox{otherwise}.
\eeqa
\end{subequations}

\begin{theorem} \label{th:ID0.6}
  Suppose that Assumption~\ref{ass:standing} holds. Then there is a
  positive number $\delta_3$ such that \eqnok{intlam} is feasible
  and bounded whenever $\delta(z,\lambda) \le \delta_3$, and its
  optimal objective is at least $\epslam$ (for $\epslam$ defined in
  \eqnok{def.dellam}). Moreover, there is a constant $\beta'>0$ such
  that $\delta(z,\hat{\lambda}) \le \beta' \delta(z,\lambda)^{\tau}$.
\end{theorem}
\begin{proof}
Let $\lambda^* \in \cSlam$ be chosen so that $\epslam =
\min_{i \in \cB_+} \lambda^*_i$. We show first that
$(\hat{t}, \hat{\lambda}) = (\epslam, \lambda^*)$ is feasible
for \eqnok{intlam}, thereby proving that
this linear program is feasible and that the optimum objective value
is at least $\epslam$.

Initially we set $\delta_3 = \delta_2$.  By Definition
\eqnok{def.dellam}, the constraint \eqnok{intlam.b} is 
satisfied by $(\hat{t}, \hat{\lambda}) = (\epslam, \lambda^*)$.  Since
$\delta(z,\lambda) \le \delta_3 = \delta_2$, we have from
Theorem~\ref{th:ID0.3} that $\cA_+ = \cB_+$, so that
\eqnok{intlam.4} also holds. Satisfaction of \eqnok{intlam.c}
follows from \eqnok{ID0.2.2}, by choice of $\delta_2$. Moreover, it is
clear from $\cA_+ = \cB_+$ that the optimal $(\hat{t}, \hat{\lambda})$
will satisfy $\hat{t} = \min_{i \in \cB_+} \hat{\lambda}_i$.

We now show that the problem  \eqnok{intlam} is bounded for 
$\delta(z,\lambda)$ sufficiently small. Let $\bar{y}$ be the vector in 
\eqnok{mfcq}, and decrease $\delta_3$ if necessary so that we can
choose a number $\zeta >0$ such that
\beq \label{eq:45}
\delta(z,\lambda) \le \delta_3  \;\; \Rightarrow \;\;
\bar{y}^T \nabla g_i (z) \le -\zeta, \;\; \mbox{for all $i \in \cA_+ = \cB_+$}.
\eeq
From the constraints \eqnok{intlam.c} and the triangle inequality, 
we have that
\beqas
\left\|  \sum_{i \in \cA_+} \hlam_i \bar{y}^T \nabla g_i(z) \right\|_1
& \le &
\| \bar{y}^T \nabla \phi(z) \|_1 + 
\left\| \bar{y}^T \nabla \phi(z) + 
\sum_{i \in \cA_+} \hlam_i \bar{y}^T \nabla g_i(z) \right\|_1 \\
& \le & \| \bar{y} \|_1 \left\| \nabla \phi(z) \right\|_{\infty} +
\| \bar{y} \|_1 
\left\| \nabla \phi(z) + \sum_{i \in \cA_+} \hlam_i \nabla g_i(z)
\right\|_{\infty} \\
& \le & \| \bar{y} \|_1 \left\| \nabla \phi(z) \right\|_{\infty} +
\| \bar{y} \|_1 \eta(z,\lambda)^{\tau}.
\eeqas
However, from \eqnok{eq:45} and $\hlam_i \ge 0$, $i \in \cA_+$, we have that
\[
\left\|  \sum_{i \in \cA_+} \hlam_i \bar{y}^T \nabla g_i(z) \right\|_1
\ge \left\| \hlam_{\cA_+} \right\|_1 \zeta.
\]
By combining these bounds, we obtain that 
\[
\left\| \hlam_{\cA_+} \right\|_1 \le \zeta^{-1} \| \bar{y} \|_1 \left[
\left\| \nabla \phi(z) \right\|_{\infty} +  \eta(z,\lambda)^{\tau} \right],
\]
whenever $\delta(z,\lambda) \le \delta_3$,
so that the feasible region for \eqnok{intlam} is bounded, as claimed.

To prove our final claim that $\delta(z,\hat{\lambda}) \le \beta'
\delta(z,\lambda)^{\tau}$ for some $\beta' >0$, we use
Theorem~\ref{th:dlamz}. We have from \eqnok{intlam.c} and the cited
theorem that
\[
\left\| \cL_z(z,\hat{\lambda}) \right\|_{\infty} \le
\eta(z,\lambda)^{\tau} \le \kappa_1^{\tau} \delta(z,\lambda)^{\tau}.
\]
For $i \in \cA_+ = \cB_+$, we have from $\hat{\lambda}_i \ge \epslam$
and $g_i(z^*) = 0$ that
\beqas
i \in \cA_+ \;\; \Rightarrow \;\;  
\left| \min(\hat{\lambda}_i, -g_i(z)) \right| 
& \le & | g_i(z)| \le |g_i(z) - g_i(z^*) | \\
&=& O(\| z-z^*\|) = O(\delta(z,\lambda)).
\eeqas
For $i \notin \cA_+$, we have $\hat{\lambda}_i=0$ and $g_i(z^*) \le 0$, and so
\beqas
i \notin \cA_+ \;\; \Rightarrow \;\;  
\left| \min(\hat{\lambda}_i, -g_i(z)) \right| 
& = & \max(0, g_i(z)) \le | g_i(z) - g_i(z^*) | \\
&=& O(\| z-z^*\|) = O(\delta(z,\lambda)).
\eeqas
By substituting the last three bounds into \eqnok{def.eta} and
applying Theorem~\ref{th:dlamz}, we obtain the result.
\end{proof}


\subsection{Computational Aspects} \label{subsec:comput}

Solution of the linear programs \eqnok{B0est} is in general less
expensive than solution of the quadratic programs or complementarity
problems that must be solved at each step of an optimization algorithm
with rapid local convergence. Linear programming software is easy to
use and readily available. Moreover, given a point $(z,\lambda)$ with
$\delta(z,\lambda)$ small, we can expect $\hat{\cA}_{\rm init}$ not
to contain many more indices than the weakly active set $\cB_0$,
so that few iterations of the ``repeat'' loop in Procedure ID0 should
be needed.

Finally, we note that when more than one iteration of the ``repeat''
loop is needed in Procedure ID0, the linear programs to be solved at
successive iterations differ only in the cost vector in
\eqnok{B0est.1}. Therefore, if the dual formulation of \eqnok{B0est}
is used, the solution of one linear program can typically be obtained
at minimal cost from the solution of the previous linear program in
the sequence.  To clarify this claim, we simplify notation and write
\eqnok{B0est} as follows:
\beq \label{lpid0.p}
\max c^T \pi \sgap \mbox{subject to} \;\; b_1 \le A \pi \le b_2, \;\; 
\pi \ge 0,
\eeq
where $\pi = [ \lambda_i ]_{i \in \cA(z,\lambda)}$, while $c$, $b_1$,
$b_2$, and $A$ are defined in obvious ways. In particular, $c$ is a
vector with elements $0$ and $1$, with the $1$'s in positions
corresponding to the index set $\hat{\cA}$. The dual of
\eqnok{lpid0.p} is
\beqas
& \max b_1^T y_1 + b_2^T y_2 \sgap \mbox{subject to} \\
& \left[ \begin{array}{ccc} A^T & -A^T & I \end{array} \right]
\bmat{c} y_1 \\ y_2 \\ s \emat = -c, \sgap
(y_1,y_2,s) \ge 0.
\eeqas
When the set $\hat{\cA}$ is changed, some of the $1$'s in the vector
$c$ are replaced by zeros. When only a few such changes are made, and
the previous optimal basis is used to hot-start the method, we expect
that only a few iterations of the dual simplex method will be needed
to recover the solution of the new linear program.

\section{SQP and Stabilized SQP} \label{sec:sqp}

In the best-known form of the SQP algorithm (with exact second-order
information), the following inequality constrained subproblem is
solved to obtain the step $\Dz$ at each iteration:
\beqa \label{sqp} 
& \min_{\Dz} \, \Dz^T \nabla  \phi(z) + \half \Dz^T
\cL_{zz} (z,\lambda) \Dz, \\ \nonumber & \makebox{\rm subject to
$g(z) + \nabg (z)^T \Dz \le 0$,} 
\eeqa 
where $(z,\lambda)$ is the current primal-dual iterate.  Denoting the
Lagrange multipliers for the constraints in \eqnok{sqp} by
$\lambda^+$, we see that the solution $\Dz$ satisfies the following
KKT conditions (cf. \eqnok{kkt.N}):
\beq \label{sqp.kkt.N}
\bmat{c} \cL_{zz} (z,\lambda) \Dz + \nabla \phi(z) + \nabg(z) \lambda^+ \\
g(z) + \nabg(z)^T \Dz 
\emat \in
\bmat{c} 0 \\ N(\lambda^+) \emat,
\eeq
where $N(\cdot)$  is defined as in \eqnok{def.Nl}.

In the stabilized SQP method, we choose a parameter $\mu \ge 0$ and
seek a solution of the following minimax subproblem for $(\Dz,
\lambda^+)$ such that $(\Dz, \lambda^+ - \lambda)$ is small:
\beqa \label{ssqp.1}
\lefteqn{ \min_{\Dz} \, \max_{\lambda^+ \ge 0} \, 
\Dz^T \nabla \phi(z) + \half \Dz^T \cL_{zz} (z,\lambda) \Dz } \\
&& \nonumber
+ (\lambda^+)^T [g(z) + \nabg(z)^T \Dz] - 
\half \mu \| \lambda^+ - \lambda \|^2.
\eeqa
The parameter $\mu$ can depend on an estimate of
the distance $\delta(z,\lambda)$ to the primal-dual solution set; for example,
 $\mu = \eta(z,\lambda)^{\sigma}$ for some $\sigma \in
(0,1)$. We can also write \eqnok{ssqp.1} as a linear complementarity 
problem, corresponding to 
\eqnok{sqp.kkt.N}, as follows:
\beq \label{ssqp.kkt.N}
\bmat{c} \cL_{zz} (z,\lambda) \Dz + \nabla \phi(z) + \nabg(z) \lambda^+ \\
g(z) + \nabg(z)^T \Dz - \mu(\lambda^+ - \lambda)
\emat \in
\bmat{c} 0 \\ N(\lambda^+) \emat.
\eeq
Li and Qi~\cite{LiQ00a} 
derive a 
quadratic program in $(\Dz, \lambda^+)$ that is equivalent to 
\eqnok{ssqp.1} and \eqnok{ssqp.kkt.N}:
\beqa
\label{ssqp.reformQP}
& \min_{(\Dz, \lambda^+)} \, \Dz^T \nabla  \phi(z) + \half \Dz^T
\cL_{zz} (z,\lambda) \Dz + \half \mu \| \lambda^+\|^2, \\ 
\nonumber & \makebox{\rm subject to $g(z) + \nabg (z)^T \Dz  - 
\mu(\lambda^+ - \lambda) \le 0$} .
\eeqa

Under conditions stronger than those assumed in this paper, the
results of Wright~\cite{Wri98a} and Hager~\cite{Hag97a} can be used to
show that the iterates generated by \eqnok{ssqp.1} (or
\eqnok{ssqp.kkt.N} or \eqnok{ssqp.reformQP}) yield superlinear
convergence of the sequence $(z^k,\lambda^k)$ of Q-order
$1+\sigma$. Our aim in the next section is to add a strategy for
adjusting the multiplier, with a view to obtaining superlinear
convergence under a weaker set of conditions.

\section{Multiplier Adjustment and Superlinear Convergence} \label{sec:adj}

We show in this section that through use of Procedure ID0 and the
multiplier adjustment strategy \eqnok{intlam}, we can devise a
stabilized SQP algorithm that converges superlinearly whenever the
initial iterate $(z^0,\lambda^0)$ is sufficiently close to the
primal-dual solution set $\cS$. Only Assumption~\ref{ass:standing} is
needed for this result.

Key to our analysis is Theorem~1 of Hager~\cite{Hag97a}. We state this
result in Appendix~\ref{sec:hager}, using our current notation and
making a slight correction to the original statement.  Here we state
an immediate corollary of Hager's result that applies under our
standing assumption.
\begin{corollary} \label{co:hager}
Suppose that Assumption~\ref{ass:standing} holds, and let $\lambda^*
\in \cSlam$ be such that $\lambda^*_i>0$ for all $i \in \cB_+$. Then
for any sufficiently large positive $\sigma_0$, there are positive
constants $\rho_0$, $\sigma_1$, $\gamma \ge 1$, and $\bar{\beta}$
such that $\sigma_0 \rho_0 < \sigma_1$, with the following property:
For any $(z^0, \lambda^0)$ with 
\beq \label{eq:60}
\|(z^0,\lambda^0) - (z^*, \lambda^*)\| \le \rho_0, 
\eeq
we can generate an iteration sequence $\{ (z^k,\lambda^k) \}$,
$k=0,1,2,\dots$, by setting
\[
(z^{k+1}, \lambda^{k+1}) = (z^k+ \Dz, \lambda^+),
\]
where, at iteration $k$, $(\Dz, \lambda^+)$ is the local solution of
the sSQP subproblem with 
\beq \label{eq:61}
(z,\lambda) = (z^k,\lambda^k), \sgap 
\mu = \mu_k \in [\sigma_0 \| z^k-z^* \|, \sigma_1],
\eeq
that satisfies
\beq \label{eq:gamma}
\left\| (z^k+ \Dz, \lambda^+) - (z^*,\lambda^*) \right\| \le \gamma
\left\| (z^0, \lambda^0) - (z^*,\lambda^*) \right\|.
\eeq
Moreover, we have
\beq \label{eq:62}
\delta(z^{k+1}, \lambda^{k+1}) \le \bar{\beta} \left[
\delta(z^k\lambda^k)^2 + \mu_k \delta(\lambda^k) \right].
\eeq
\end{corollary}

Recalling our definition \eqnok{def.dellam} of $\epslam$, we define
the following parametrized subset of $\cSlam$:
\beq \label{def.nuSlam}
\cSlam^{\nu} \defeq 
\{ \lambda \in \cSlam \, | \,
\min_{i \in \cB_+} \lambda_i \ge \nu \epslam \}.
\eeq
It follows easily from the MFCQ assumption and \eqnok{def.dellam} that
$\cSlam^{\nu}$ is nonempty, closed, bounded, and therefore compact for
any $\nu \in [0,1]$.

We now show that the particular choice of stabilization parameter $\mu
= \eta(z,\lambda)^\sigma$, for some $\sigma \in (0,1)$, eventually
satisfies \eqnok{eq:61}.
\begin{lemma} \label{lem:choosemu}
Suppose the assumptions of Corollary~\ref{co:hager} are satisfied, and
let $\lambda^*$ be as defined there. Let $\sigma$ be any constant in $(0,1)$.
Then there is a quantity
$\rho_2 \in (0,\rho_0]$ such that when $(z^0,\lambda^0)$ satisfies
\beq \label{eq:64}
\|(z^0,\lambda^0) - (z^*, \lambda^*)\| \le \rho_2, 
\eeq
the results of Corollary~\ref{co:hager} hold when we set the
stabilization parameter at iteration $k$ to the following particular
value:
\beq \label{eq:65}
\mu = \mu_k = \eta(z^k,\lambda^k)^{\sigma}.
\eeq
\end{lemma}
\begin{proof}
We prove the result by showing that $\mu_k$ defined by \eqnok{eq:65}
satisfies \eqnok{eq:61} for some choice of $\rho_2$. For
contradiction, suppose that no such choice of $\rho_2$ is possible,
so that for each $\ell=1,2,3,\dots$, there is a starting point
$(\zsl^0, \lsl^0)$ with
\beq \label{eq:67}
\left\| (\zsl^0, \lsl^0) - (z^*,\lambda^*) \right\| \le \ell^{-1} \rho_0
\eeq
such that the sequence $\left\{ \left( \zsl^k, \lsl^k \right)
\right\}_{k=0,1,2,\dots}$ generated from this starting point in the
manner prescribed by Corollary~\ref{co:hager} with $\mu_k =
\eta(\zsl^k, \lsl^k)^{\sigma}$ eventually comes across an index
$k_{\ell}$ such that this choice of $\mu_k$ violates \eqnok{eq:61},
that is, one of the following two conditions holds:
\begin{subequations} \label{eq:68}
\beqa \label{eq:68a}
\sigma_0 \left\|  \zsl^{k_{\ell}} - z^* \right\| 
& > & \eta(\zsl^k, \lsl^k)^{\sigma}, \\
\label{eq:68b}
\sigma_1 & < & \eta(\zsl^k, \lsl^k)^{\sigma}.
\eeqa
\end{subequations}
Assume that $k_{\ell}$ is the first such index for which the violation
\eqnok{eq:68} occurs. By \eqnok{eq:gamma} and \eqnok{eq:67}, we have
that
\beq \label{eq:69}
\left\| \left( \zsl^{k_{\ell}}, \lsl^{k_{\ell}} \right) 
- (z^*,\lambda^*) \right\| \le 
\gamma \left\| \left(  \zsl^0, \lsl^0 \right) - (z^*,\lambda^*) \right\| 
\le \gamma \ell^{-1} \rho_0.
\eeq
Therefore by Theorem~\ref{th:dlamz} and \eqnok{errdef}, we have for
$\ell$ sufficiently large that
\beqa
\nonumber
\frac{\eta \left( \zsl^{k_{\ell}}, \lsl^{k_{\ell}} \right)^{\sigma}}{
\left\| \zsl^{k_{\ell}} - z^* \right\|} & \ge & 
\frac{\eta \left( \zsl^{k_{\ell}}, \lsl^{k_{\ell}} \right)^{\sigma}}{
\delta \left( \zsl^{k_{\ell}}, \lsl^{k_{\ell}} \right)} \\
\nonumber
& \ge &
\kappa_0^{\sigma} 
\delta \left( \zsl^{k_{\ell}}, \lsl^{k_{\ell}} \right)^{\sigma-1} \\
\nonumber
& \ge &
\kappa_0^{\sigma} 
\left\| \left( \zsl^{k_{\ell}}, \lsl^{k_{\ell}} \right) - (z^*,\lambda^*)
\right\|^{\sigma-1} \\
\label{eq:69b}
& \ge &
\kappa_0^{\sigma} \gamma^{\sigma-1} \rho_0^{\sigma-1} \ell^{1-\sigma}.
\eeqa
Hence, taking limits as $\ell \uparrow \infty$, we have that 
\[
\frac{\eta \left( \zsl^{k_{\ell}}, \lsl^{k_{\ell}} \right)^{\sigma}}{
\left\| \zsl^{k_{\ell}} - z^* \right\|}  \to \infty \sgap
\mbox{as $\ell \uparrow \infty$.}
\]
Dividing both sides of \eqnok{eq:68a} by $\left\| \zsl^{k_{\ell}} -
z^* \right\|$, we conclude from finiteness of $\sigma_0$ that
\eqnok{eq:68a} is impossible.

By using Theorem~\ref{th:dlamz} again together with \eqnok{eq:69}, we
obtain
\beqas
\eta \left( \zsl^{k_{\ell}}, \lsl^{k_{\ell}} \right) & \le &
\kappa_1 \delta \left( \zsl^{k_{\ell}}, \lsl^{k_{\ell}} \right) \\
& \le & \kappa_1 \left\| 
\left( \zsl^{k_{\ell}}, \lsl^{k_{\ell}} \right) - (z^*,\lambda^*)
\right\| \\
& \le &
\kappa_1 \gamma \rho_0 \ell^{-1},
\eeqas
and therefore $\eta \left( \zsl^{k_{\ell}}, \lsl^{k_{\ell}}
\right)^{\sigma} \to 0$ as $\ell \uparrow \infty$. Hence,
\eqnok{eq:68b} cannot occur either, and the proof is complete.
\end{proof}

We now use a compactness argument to extend Corollary~\ref{co:hager}
from the single multiplier $\lambda^*$ in the relative interior of
$\cSlam$ to the entire set $\cSlam^{\nu}$, for any $\nu \in (0,1]$.
\begin{theorem} \label{th:ssqp2}
  Suppose that Assumption~\ref{ass:standing} holds, and fix $\nu \in
  (0,1]$. Then there are positive constants $\hat{\delta}$, $\gamma \ge 1$,
  and $\beta$ such that the following property holds: Given
  $(z^0,\lambda^0)$ with
\[
\makebox{\rm dist} \left( (z^0,\lambda^0), \cSlam^{\nu} \right) 
\le \hat{\delta}, 
\]
the iteration sequence $\{ (z^k,\lambda^k) \}_{k=0,1,2,\dots}$
generated in the manner described in Corollary~\ref{co:hager}, with
$\mu_k$, $k=0,1,2\dots$ chosen according to \eqnok{eq:65}, satisfies
the following relations:
\begin{subequations}
\beqa
\label{eq:70}
\delta(z^{k+1},\lambda^{k+1}) & \le & \beta \delta(z^k,\lambda^k)^{1+\sigma} \\
\label{eq:71}
\lambda^k_i & \ge & \frac12 \nu \epslam, \sgap
\mbox{for all $i \in \cB_+$ and all $k=0,1,2\dots$}.
\eeqa
\end{subequations}
\end{theorem}
\begin{proof}
For each $\lambda^* \in \cSlam^{\nu}$, we use Corollary~\ref{co:hager}
to obtain positive constants $\sigma_0(\lambda^*)$ (sufficiently
large), $\sigma_1(\lambda^*)$, $\gamma(\lambda^*)$, and $\bar{\beta}
(\lambda^*)$, using the argument $\lambda^*$ for each constant to
emphasize the dependence on the choice of multiplier $\lambda^*$. In
the same vein, let $\rho_2(\lambda^*) \in (0, \rho_0(\lambda^*)]$
be the constant from Lemma~\ref{lem:choosemu}. Now choose
$\hat{\delta}(\lambda^*)>0$ for each $\lambda^* \in \cSlam^{\nu}$ 
in such a way that
\begin{subequations} \label{eq:72}
\beqa
\label{eq:72a}
& 0 < \hat{\delta} (\lambda^*) \le \half \rho_2(\lambda^*), \\
\label{eq:72b}
& \gamma(\lambda^*) \hat{\delta}(\lambda^*) \le \quarter \nu \epslam,
\eeqa
\end{subequations}
and consider the following open cover of $\cSlam^{\nu}$:
\beq \label{eq:73}
\cup_{\lambda^* \in \cSlam^{\nu}} 
\left\{ \lambda \, | \, \| 
\lambda - \lambda^* \| < \hat{\delta}(\lambda^*) \right\}.
\eeq
By compactness of $\cSlam^{\nu}$, we can find a finite subcover
defined by points $\hlam^1, \hlam^2, \dots, \hlam^f \in \cSlam^{\nu}$ 
as follows:
\beq \label{eq:74}
\cSlam^{\nu} \subset \cV \defeq \cup_{j=1,2,\dots,f} 
\left\{ \lambda \, | \, \| 
\lambda - \hlam^j \| < \hat{\delta}(\hlam^j) \right\}.
\eeq
$\cV$ is an open neighborhood of $\cSlam^{\nu}$. Now define 
\beq \label{eq:75}
\gamma \defeq \max_{j=1,2,\dots,f} \gamma(\hlam^j), \sgap
\bar{\beta} \defeq  \max_{j=1,2,\dots,f} \bar{\beta} (\hlam^j), \sgap
\delta \defeq \max_{j=1,2,\dots,f} \hat{\delta}(\hlam^j).
\eeq
Also, choose a quantity $\hat{\delta}>0$ with the following properties:
\begin{subequations} \label{eq:76}
\beqa
\label{eq:76a}
\hat{\delta} & \le & \min_{j=1,2,\dots,f} \hat{\delta} (\hlam^j) \le \delta, \\
\label{eq:76b}
\left\{ \lambda \, | \, 
\mbox{dist}(\lambda, \cSlam^{\nu}) \le \hat{\delta} \right\} & \subset & \cV, \\
\label{eq:76c}
\hat{\delta} & \le & \frac{\nu \epslam}{4 \gamma}, \\
\label{eq:76d}
\hat{\delta} & \le & 1.
\eeqa
\end{subequations}

Now consider $(z^0,\lambda^0)$ with 
\beq \label{eq:79}
\left\| (z^0,\lambda^0) - (z^*, \lambda^*) \right\| \le \hat{\delta}, \sgap
\mbox{for some $\lambda^* \in \cSlam^{\nu}$}.
\eeq
We have $\mbox{dist}(\lambda^0, \cSlam^{\nu}) \le \hat{\delta}$, and so
$\lambda^0 \in \cV$. It follows that for some $j=1,2,\dots,f$, we have
\beq \label{eq:78}
\| \lambda^0 - \hlam^j \| \le \hat{\delta} (\hlam^j).
\eeq
Moreover, since 
$\| z^0-z^* \| \le \hat{\delta}$,
we have from \eqnok{eq:76a} that
\beq \label{eq:80}
\left\| (z^0,\lambda^0) - (z^*,\hlam^j) \right\| \le \hat{\delta} +
\hat{\delta} (\hlam^j) \le 2 \hat{\delta} (\hlam^j) \le \rho_2(\hlam^j),
\eeq
where the final inequality follows from  \eqnok{eq:72a}.
Application of Corollary~\ref{co:hager} and Lemma~\ref{lem:choosemu}
now ensures that the stabilized SQP sequence starting at
$(z^0,\lambda^0)$ with $\mu=\mu_k$ chosen according to \eqnok{eq:65}
yields a sequence $\{ (z^k\lambda^k) \}_{k=0,1,2,\dots}$ satisfying
\beqa
\nonumber
\left\| (z^k,\lambda^k) - (z^*,\hlam^j) \right\| & \le &
\gamma(\hlam^j) \left\| (z^0,\lambda^0) - (z^*,\hlam^j) \right\| \\
\label{wholeseq.dist}
& \le & 2 \gamma (\hlam^j) \hat{\delta} (\hlam^j) 
\le 2 \gamma \delta,
\eeqa
where we used \eqnok{eq:75} to obtain the final inequality. 

To prove \eqnok{eq:70}, we have from Lemma~\ref{lem:choosemu},
Corollary~\ref{co:hager}, the bound \eqnok{deltas},
Theorem~\ref{th:dlamz}, the definition \eqnok{eq:75}, and the
stabilizing parameter choice \eqnok{eq:65} that
\beqas
\delta(z^{k+1}, \lambda^{k+1}) & \le &
\bar{\beta} (\hlam^j) 
\left[ \delta(z^k,\lambda^k)^2 + \mu_k \delta(\lambda^k) \right] \\
& \le &
\bar{\beta} 
\left[ \delta(z^k,\lambda^k)^2 + \eta(z^k\lambda^k)^{\sigma} 
\delta(z^k, \lambda^k) \right] \sgap \mbox{from (\protect\ref{eq:75}) and
(\protect\ref{eq:65})} \\
& \le &
\bar{\beta} 
\left[ \delta(z^k,\lambda^k)^2 + 
\kappa_1^{\sigma} \delta(z^k,\lambda^k)^{1+\sigma} \right] \sgap
\mbox{from Theorem~\protect{\ref{th:dlamz}}} \\
& \le & 
\bar{\beta} \left(  (2\gamma \delta)^{1-\sigma} +\kappa_1^{\sigma} \right) 
\delta(z^k,\lambda^k)^{1+\sigma}, 
\eeqas
where in the last line we use $\delta(z^k,\lambda^k) \le
\mbox{dist}((z^k,\lambda^k), \cSlam^{\nu}) \le 2\gamma \delta$.
Therefore, the result \eqnok{eq:70} follows by setting $\beta =
\bar{\beta} \left( (2\gamma \delta)^{1-\sigma} + \kappa_1^{\sigma}
\right)$.

Finally, we have from \eqnok{eq:72b} (with $\lambda^* = \hlam^j$) 
and \eqnok{wholeseq.dist} that
\[
\mbox{dist} \left( (z^k,\lambda^k), \cSlam^{\nu} \right) 
\le 2 \gamma (\hlam^j)  \hat{\delta} (\hlam^j) \le \frac12 \nu \epslam.
\]
Therefore, we have
\[
i \in \cB_+ \; \Rightarrow \;
\lambda^k_i \ge \min_{\lambda^* \in \cSlam^{\nu}} \lambda^*_i -
\frac12 \nu \epslam \ge
\nu \epslam - \frac12 \nu \epslam = \frac12 \nu \epslam,
\]
verifying \eqnok{eq:71} and completing the proof.
\end{proof}


We are now ready to state a stabilized SQP algorithm, in which
multiplier adjustment steps (consisting of Procedure ID0 followed by
solution of \eqnok{intlam}) are applied when the convergence does
not appear to be rapid enough.

\medskip

\btab
\> {\bf Algorithm sSQPa} \\
\> given $\sigma \in (0,1)$, $\tau$ and $\hat{\tau}$ with 
$0 < \hat{\tau} < \tau < 1$, tolerance {\tt tol}; \\
\> given initial point $(z^0,\lambda^0)$ with $\lambda^0 \ge 0$; \\
\> $k \leftarrow 0$; \\
\> calculate $\cA(z^0,\lambda^0)$ from \eqnok{def.Best}; \\
\> call Procedure ID0 to obtain $\cA_+$, $\cA_0$; 
solve \eqnok{intlam} to obtain $\hlam^0$; \\
\> $\lambda^0 \leftarrow \hlam^0$; \\
\> {\bf repeat} \\
\>\> solve \eqnok{ssqp.1} with $(z,\lambda) = (z^k,\lambda^k)$ and 
$\mu = \mu_k = \eta(z^k,\lambda^k)^{\sigma}$ \\
\>\>\> to obtain  $(\Delta z, \lambda^+)$; \\
\>\> {\bf if} $\eta(z^k + \Delta z, \lambda^+) \le 
\eta(z^k,\lambda^k)^{1+\sigma/2}$ \\
\>\>\> $(z^{k+1}, \lambda^{k+1}) \leftarrow (z^k + \Delta z, \lambda^+)$; \\
\>\>\> $k \leftarrow k+1$; \\
\>\> {\bf else} \\
\>\>\> calculate $\cA(z^k,\lambda^k)$ from \eqnok{def.Best}; \\
\>\>\> call Procedure ID0 to obtain $\cA_+$, $\cA_0$;
 solve \eqnok{intlam} to obtain $\hlam^k$; \\
\>\>\> $\lambda^k \leftarrow \hlam^k$; \\
\>\> {\bf end (if)} \\
\> {\bf until} $\eta(z^k,\lambda^k) < {\tt tol}$.
\etab

\medskip

The following result shows that when $(z^0,\lambda^0)$ is close enough
to $\cS$, the initial call to Procedure ID0 is the only one needed.

\begin{theorem} \label{th:ssqpa}
  Suppose that Assumption~\ref{ass:standing} holds. Then there is a
  constant $\bar{\delta}>0$ such that for any $(z^0,\lambda^0)$ with
  $\delta(z^0,\lambda^0) \le \bar{\delta}$, the ``if'' condition in
  Algorithm sSQPa is always satisfied, and the sequence
  $\delta(z^k,\lambda^k)$ converges superlinearly to zero with Q-order
  $1+\sigma$.
\end{theorem}
\begin{proof}
  Our result follows from Theorems~\ref{th:ID0.6} and \ref{th:ssqp2}.
  Choose $\nu = 1/2$ in Theorem~\ref{th:ssqp2}, and let
  $\hat{\delta}$, $\gamma$, and $\beta$ be as defined there. Using
  also $\delta_3$ and $\beta'$ from Theorem~\ref{th:ID0.6} and
  $\epslam$ defined in \eqnok{def.dellam}, we choose $\bar{\delta}$ as
  follows:
\beq \label{choose.bardel}
\bar{\delta} = \min \left(
\delta_3, 
\hat{\delta},
\left( \frac{\epslam}{2 \beta'} \right)^{1/\tau}, 
\left( \frac{\hat{\delta}}{\beta'} \right)^{1/\tau}, 
\frac{1}{(2 \beta)^{1/\sigma}}, 
\kappa_0 \left( \frac{\kappa_0}{\beta \kappa_1} \right)^{2/\sigma}
\right).
\eeq
Now let $(z^0,\lambda^0)$ satisfy $\delta(z^0,\lambda^0) \le
\bar{\delta}$, and let $\hat{\lambda}^0$ be calculated from
\eqnok{intlam}. From Theorem~\ref{th:ID0.6} and
\eqnok{choose.bardel}, we have that
\beq \label{eq:101}
\delta(z^0,\hlam^0) \le \beta' \delta(z^0,\lambda^0)^{\tau} 
\le \beta' \bar{\delta}^{\tau} \le \frac12 \epslam
\eeq
and
\begin{subequations} \label{eq:102}
\beqa \label{eq:102.a}
\hlam^0_i & \ge & \epslam, \sgap \mbox{for all $i \in \cB_+$}, \\
\label{eq:102.b}
\hlam^0_i &=& 0, \sgap \mbox{for all $i \notin \cB_+$}.
\eeqa
\end{subequations}
Since $\cSlam$ is closed, there is a vector $\hlam^* \in \cSlam$ such that 
\beq \label{eq:103}
\delta(z^0,\hlam^0) = \left\| (z^0,\hlam^0) - (z^*, \hlam^*) \right\|.
\eeq
From \eqnok{eq:101} and \eqnok{eq:102.a}, we have that 
\[
i \in \cB_+ \; \Rightarrow \; \hlam^*_i \ge \hlam^0_i - \frac12 \epslam \ge 
\frac12 \epslam,
\]
so that $\hlam^* \in \cSlam^{\nu}$ for $\nu = 1/2$.
We therefore have from \eqnok{eq:101},
\eqnok{eq:103}, and \eqnok{choose.bardel} that
\beq \label{eq:104}
\mbox{dist} ((z^0,\hlam^0), \cSlam^{\nu}) =
\left\| (z^0,\hlam^0) - (z^*, \hlam^*) \right\| \le \beta' \bar{\delta}^{\tau}
\le \hat{\delta}.
\eeq
From here on, we set $\lambda^0 \leftarrow \hlam^0$, as in Algorithm
sSQPa.  Because of the last bound, we can apply Theorem~\ref{th:ssqp2}
to $(z^0,\lambda^0)$. We use this result to prove the following
claims. First,
\beq \label{eq:105}
\bar{\delta} \ge \delta(z^0, \lambda^0) \ge 2 \delta(z^1, \lambda^1) \ge
4 \delta(z^2, \lambda^2) \ge \cdots.
\eeq
Second, 
\beq \label{eq:106}
\eta(z^{k+1}, \lambda^{k+1}) \le \eta(z^k, \lambda^k)^{1+\sigma/2},
\sgap \mbox{for all $k=0,1,2, \dots$}.
\eeq
We prove both claims by induction. For $k=0$ in \eqnok{eq:105}, we
have from \eqnok{eq:104} and $\bar{\delta} \le \hat{\delta}$ in
\eqnok{choose.bardel} that $\delta(z^0,\lambda^0) \le \bar{\delta}$.
Assume that the first $k+1$ inequalities in
\eqnok{eq:105} have been verified.
From \eqnok{eq:70} and \eqnok{choose.bardel}, we have that
\[
\delta(z^{k+1}, \lambda^{k+1}) \le \beta
\delta(z^k,\lambda^k)^{1+\sigma} \le \beta \bar{\delta}^{\sigma}
\delta(z^k,\lambda^k) \le \frac12 \delta(z^k,\lambda^k),
\]
so that the next inequality in the chain is also satisfied. 
For \eqnok{eq:106}, we have from
Theorem~\ref{th:dlamz}, \eqnok{eq:70}, and \eqnok{eq:105} that
\beqas
\eta(z^{k+1}, \lambda^{k+1}) & \le & \kappa_1 \delta(z^{k+1}, \lambda^{k+1}) \\
& \le & \beta \kappa_1 \delta(z^k,\lambda^k)^{1+\sigma} \\
& \le & \beta \kappa_1 \bar{\delta}^{\sigma/2} 
\delta(z^k,\lambda^k)^{1+\sigma/2} \\
& \le & \beta \kappa_1 \bar{\delta}^{\sigma/2}  \kappa_0^{-1-\sigma/2} 
\eta (z^k,\lambda^k)^{1+\sigma/2} \\
& \le & \eta (z^k,\lambda^k)^{1+\sigma/2},
\eeqas
where the last bound follows from \eqnok{choose.bardel}.  Hence,
\eqnok{eq:106} is verified, so that the condition in the ``if''
statement of Algorithm sSQPa is satisfied for all $k=0,1,2,\dots$.
Superlinear convergence with Q-order $1+\sigma$ follows from
\eqnok{eq:70}.
\end{proof}

\section{Summary and Possible Extensions} \label{sec:ext}

We have presented a technique for identifying the active inequality
constraints at a local solution of a nonlinear programming problem,
where the standard assumptions---existence of a strictly complementary
solution and linear independence of active constraints gradients---are
replaced by weaker assumptions. We have embedded this technique in a
stabilized SQP algorithm, resulting in a method that converges
superlinearly under the weaker assumptions when started at a point
sufficiently close to the (primal-dual) optimal set. 

The primal-dual algorithm described by Vicente and
Wright~\cite{VicW00} can also be improved by using the techniques
outlined here. In that paper, strict complementarity is assumed along
with MFCQ, and superlinear convergence is proved provided both
$\delta(z^0,\lambda^0)$ is sufficiently small and $\lambda^0_i \ge
\gamma$, for all $i \in \cB = \cB_+$ and some $\gamma>0$. If we apply
the active constraint detection procedure \eqnok{def.Best} and the
subproblem \eqnok{intlam} to {\em any} initial point $(z^0,\lambda^0)$
with $\delta(z^0,\lambda^0)$ sufficiently small, the same convergence
result can be obtained without making the positivity assumption on the
components of $\lambda^0_{\cB_+}$. (Because of the strict
complementarity assumption, Procedure ID0 serves only to verify that
$\cB = \cB_+$.)

Numerous issues remain to be investigated. We believe that degeneracy
is an important issue, given the large size of many modern
applications of nonlinear programming and their nature as
discretizations of continuous problems. Nevertheless, the practical
usefulness of constraint identification and stabilization techniques
remains to be investigated. The numerical implications should also be
investigated, since implementation of these techniques may require
solution of ill-conditioned systems of linear equations (see
M.~H.~Wright~\cite{MWri98a} and S.~J.~Wright~\cite{Wri98b}). Embedding
of these techniques into globally convergence algorithmic frameworks
needs to be examined. We should investigate generalization to equality
constraints, possibly involving the use of the ``weak'' MFCQ
condition, which does not require linear independence of the equality
constraint gradients.

%
%

\section*{Acknowledgments}

We thank Bill Hager for discussions of his key result,
Theorem~\ref{th:hager}.

\appendix

\section{Hager's Theorem} \label{sec:hager}

We restate Theorem~1 of Hager~\cite{Hag97a}, making a slight correction to
the original statement concerning the conditions on
$(z^0,\lambda^0)$ and the radius of the neighborhood containing the
sequence $\{ (z^k,\lambda^k) \}$. No modification to Hager's analysis
is needed to prove the following version of this result.
\begin{theorem} \label{th:hager}
Suppose that $z^*$ is a local solution of \eqnok{nlp}, and that $\phi$ and
  $g$ are twice Lipschitz continuously differentiable in a
  neighborhood of $z^*$. Let $\lambda^*$ be some multiplier such that the
  KKT conditions \eqnok{kkt} are satisfied, and define
\[
\bar{\cB} \defeq \{ i \, | \, \lambda^*_i >0 \}.
\]
Suppose that there is an $\alpha>0$ such that 
\[
w^T \cL_{zz}(z^*, \lambda^*)w \ge \alpha \| w \|^2, \;\; 
\mbox{for all $w$ such that $\nabla g_i(z^*)^Tw=0$, for all $i \in \bar{\cB}$}.
\]
Then for any choice of $\sigma_0$ sufficiently large, there are
positive constants $\rho_0$, $\sigma_1$, $\gamma \ge 1$, and
$\bar{\beta}$ such that $\sigma_0 \rho_0 < \sigma_1$, with the
following property: For any $(z^0, \lambda^0)$ with 
\[
\|(z^0,\lambda^0) - (z^*, \lambda^*)\| \le \rho_0, 
\]
we can generate an iteration sequence $\{ (z^k,\lambda^k) \}$,
$k=0,1,2,\dots$, by setting
\[
(z^{k+1}, \lambda^{k+1}) = (z^k+ \Dz, \lambda^+),
\]
where, at iteration $k$, $(\Dz, \lambda^+)$ is the local solution of
the sSQP subproblem with 
\[
(z,\lambda) = (z^k,\lambda^k), \sgap 
\mu = \mu_k \in [\sigma_0 \| z^k-z^* \|, \sigma_1],
\]
that satisfies
\[
\left\| (z^k+ \Dz, \lambda^+) - (z^*,\lambda^*) \right\| \le \gamma
\left\| (z^0, \lambda^0) - (z^*,\lambda^*) \right\|.
\]
Moreover, we have
\[
\delta(z^{k+1}, \lambda^{k+1}) \le \bar{\beta} \left[
\delta(z^k\lambda^k)^2 + \mu_k \delta(\lambda^k) \right].
\]
\end{theorem}

\bibliography{refs}
\bibliographystyle{siam}

\end{document}